\noindent{\bf On the denseness of the invertible group in
Banach algebras}

\vskip 20pt
\centerline{T. W. Dawson and J. F. Feinstein}
\vskip 20pt


\font\bbroman=msbm10
\font\smallbbroman=msbm7
\font\tinybbroman=msbm5
\textfont5=\bbroman \scriptfont5=\smallbbroman
 \scriptscriptfont5=\tinybbroman
\def\blackboardroman{\fam5 \bbroman}

 \def\C{{\blackboardroman C}}	
 \def\N{{\blackboardroman N}}	
 \def\Z{{\blackboardroman Z}}	
 \def\R{{\blackboardroman R}}	
 \def\F{{\blackboardroman F}}	

\def\sqr#1#2{{\vcenter{\hrule height.#2pt
	\hbox{\vrule width.#2pt height#1pt \kern #1pt
		\vrule width.#2pt}
	\hrule height.#2pt}}}
\def\square{\mathchoice\sqr77\sqr77\sqr77\sqr77}

\def\eop{\hfill $\square$ \vskip10pt}
\def\thm{\vskip 10pt\noindent THEOREM }
\def\dfn{\vskip 10pt\noindent DEFINITION }
\def\ex{\vskip 10pt\noindent EXAMPLE }
\def\cor{\vskip 10pt\noindent COROLLARY }

\def\prop{\vskip 10pt\noindent PROPOSITION }
\def\npar{\vskip 10pt\noindent}
\def\pf{\par\noindent {\it Proof.} }
\def\eop{$\hfill\square$}
\def\ref#1{[{\bf #1}]}

\def\pri{{^\prime}}

\def\inv{{^{-1}}}

\def\hA{\hat A}			
\def\O{\Omega}			
\def\o{\omega}
\def\a{\alpha}

\def\norm#1{\left\Vert#1\right\Vert}
\def\ab#1{\left\vert#1\right\vert}

\def\st{\,\colon\;}
\def\set#1{\left\{#1\right\} }

\def\a{\alpha}
\def\Aa{A_{\alpha}}

\def\e{\varepsilon}

\def\o{\omega}
\def\d{\delta}

\def\tsr{\hbox{tsr}\,}

\def\xb{\bar x}

\def\sp#1{\sigma\left(#1\right)}

\def\int{^\circ}


\def\dim{{\rm dim}\,}
\def\tsr{{\rm tsr}\,}
\def\bGA{\overline{G(A)}}
\def\d{\delta}

\def\res{{\rm res}\,}

\def\tc#1{\tilde c_{n-#1}}
\def\tcb#1{\tilde c_{n-(#1)}}
\def\prin{\bigl(\a(x)\bigr)}

\def\b{\beta}

\def\bpr{b^\prime}
\def\tA{\widetilde{A}}
\def\bppr{b^{\prime\prime}}

\def\F{{\cal F}}



\vskip 15pt
\centerline{ABSTRACT}
\vskip 7pt
{\sl\noindent We examine the condition that a complex Banach algebra
$A$ have dense invertible group. We show that,
for commutative algebras, this property is preserved by
integral extensions. We also investigate the connections with an old problem
in the theory of uniform algebras.\ \rm [MSC: 46J10, 46H05]}


\vskip 15pt
\noindent{\bf Introduction}
\vskip 10pt

\noindent In this paper we investigate the condition that a complex Banach algebra, $A$, have a
dense invertible group, $G(A)$. Throughout we shall consider only unital, complex algebras;
usually they will also be commutative.

The condition was first explored by Robertson in $\ref{22}$ for
$C^*$-algebras.
As remarked there,
a complete characterisation of it in terms of the topology of the
maximal ideal space has been known for
some time in the commutative case. It is so fundamental that we state it now: if
$X$ is a compact Hausdorff space then $C(X)$ has dense invertible group if and only
if $\dim X$, the covering dimension of $X$, is not more than $1$. For a proof
of this fact see for example $\ref{20}$, Proposition 3.3.2.
This book is a useful source of information
on dimension.

The condition has since attracted considerable
(implicit) attention as the first case of the condition `$\tsr A\le m$' where $m$ is
a natural number, $A$ a topological algebra, and $\tsr A$ its
`topological stable rank'. This notion was introduced in the seminal
paper $\ref {21}$ and for a very useful survey and collection of results
on this subject we refer the reader to $\ref 5$.

We also mention that the condition has been used to characterise
commutative $C^*$-algebras in which all elements have
square roots ($\ref {13}$). Furthermore, in Chapter 10 of $\ref{20}$,
the
dimension of a compact Hausdorff space, $X$, is described precisely in terms of
the set of (uniformly closed) subalgebras of $C(X,\R)$
which are integrally closed in $C(X,\R)$.
(Recall that an element $b$ of an algebra, $B$, is said to be {\it integral} over
a unital subalgebra, $A$, if there exists a monic polynomial with coefficients in $A$
which has $b$ as a root.) However we shall only consider complex
Banach algebras here.

In general it is not easy to determine the closure of the group of invertible elements
of a Banach algebra. Again we quote an example from $\ref{22}$: when $A$ is the disc algebra,
$\bGA$ is equal to the set of elements of $A$ whose zero set lies in the unit circle.
Further introductory examples are given in Section 1 of this paper.

Our study of the denseness of the invertible group
was prompted by an observation (Theorem 1.7) which sheds light on the
old open problem of whether or not there is a non-trivial uniform algebra on
the unit interval whose maximal ideal space is equal to the unit interval.
However our main result, proved in the section after, is
that the property of having a dense invertible
group is inherited by integral extensions of a commutative Banach algebra.


\vskip 12pt
\noindent{\bf 1. Notation, Preliminary Examples and Results}
\vskip 10pt

\noindent We begin by setting out some notation and standard terminology.

The Gelfand transform of a normed algebra, $A$, will be denoted $\hA$. The algebra $A$
is said to be {\it symmetric} if $\bar f\in \hA$ whenever $f\in\hA$.
We shall use the notation $B_A(a,r)$ for the open ball of centre $a$
and radius $r$ in the normed space $A$. The unit circle is denoted $S^1$; we regard
it as a subset of the complex plane.

\npar The following
definition, which seems quite standard in the literature, is useful.

\dfn 1.1. A subalgebra, $B$, of an algebra,
$A$, is called {\it full} if we have $G(B)=B\cap G(A)$.

\npar For example, $R_0(X)$, the algebra of rational functions on a compact plane set, $X$, with
poles off $X$ and $C^n(I)$ are both full subalgebras of their completions
with respect to the uniform norm.

\ex 1.2. Let $B$ be the restriction to $S^1$ of the
algebra of complex rational functions with poles off $S^1\cup\set2$.
Then $B$ is easily seen to be
uniformly dense in $A=C(S^1)$ (by the
Stone-Weierstrass theorem for example). However, $B$ is not a full subalgebra of $A$ because, for
example, the function given by $z-2\in B$ is invertible in $A$. \eop

\npar For future reference we note the following elementary result.

\prop 1.3. Let $A$ be a Banach algebra with maximal ideal space $\O$
such that $\hA$
is dense in $C(\O)$ (for example if $A$ is symmetric). Give $\hA$ the uniform
norm. Then $\hA$ has dense invertible
group if and only if $C(\O)$ has.

\pf It is elementary that $\hA$ is a full subalgebra of $C(\O)$. Moreover, if
$\hA$ is dense in $C(\O)$ then, since $G(C(\O))$ is open in $C(\O)$, we have from
elementary topology that
$$\overline{G(\hA)}=\overline{\hA\cap G(C(\O))}=\overline{G(C(\O))}$$
where the closures are taken in $C(\O)$. This shows that if $G(C(\O))$ is dense in
$C(\O)$ then the invertible elements of $\hA$ are dense in $\hA$. The converse
is trivial.
\eop

\npar It follows from facts mentioned in Section 1.2 of $\ref5$ that if $A$
is a commutative Banach algebra with maximal
ideal space $\O$ and $A$ is regular (that is, the hull-kernel topology on $\O$ is the
Gelfand topology) then the denseness of $G(A)$ in $A$ implies that $C(\O)$ also
has dense invertible group. We do not know if the converse is true for
regular Banach algebras. Neither do we know any examples of non-regular Banach
algebras for which $G(A)$ is dense in $A$ but $C(\O)$ does not have dense invertible
group.

\ex 1.4. In $\ref{21}$ it is shown that if $A$ is a unital $C^*$-algebra then $A$ has
dense invertible group if and only if $M_n(A)$ has for every $n\in\N$. We remark that
the proof given in $\ref{21}$ applies to all Banach algebras.

In view of the fact that $M_n(A)=M_n(\C)\hat\otimes A$ (see for example $\ref{19}$, p. 43), where
$\hat\otimes$ is the projective tensor product, one
might be tempted to conjecture that if $A_1$ and $A_2$ were commutative unital Banach
algebras with dense invertible groups then
$A_1\hat\otimes A_2$ would have dense invertible
group. However if we let $A_1=A_2=C(I)$ then by a standard result (see for
example $\ref 6$, Proposition 2.3.7) their tensor product, $B$, has a
maximal ideal space homeomorphic to $I\times I$ (dimension 2) and $B$ is easily seen
to be symmetric. Therefore $B$ could not have dense invertible
group by Proposition 1.3 above.

\prop 1.5. Suppose that $A$ is a commutative unital Banach algebra which is
rationally generated by a single element, $a$. (See $\ref{6}$, Definition 2.2.7:
we mean that the algebra
of elements obtained
by applying the rational functions with poles off $\sp a$ to $a$ is dense in $A$.) If
$\sp a$ has empty interior then $A$ has dense invertible group.

\pf This follows from the fact that if $K$ is a subset of $\C$ with
empty interior and $\phi$ is a function holomorphic on a neighbourhood
of $K$ then $\phi(K)$ also has no interior. Thus $A$ has a dense set
of elements whose spectra have no interior, which, as is noted in
$\ref{7}$, is equivalent to the condition $\bGA=A$.\eop

\npar For the definition of a weight sequence and Beurling algebras, which appear
in the next example, we refer the reader to $\ref6$, Definition 4.6.6.

\cor 1.6. Ler $\o$ be a weight sequence on $\Z$ and let
$A$ be the associated Beurling algebra, $\ell^1(\Z,\o)$. Then $A$ has dense invertible
group if and only if $\rho_-:=\lim_{n\to+\infty}\o_n^{1/n}=\lim_{n\to+\infty}\o_{-n}^{1/{-n}}
=:\rho_+$.
\pf It is standard (see for example $\ref {6}$, Theorem 4.6.7) that $A$ is
unitally polynomially generated by $\set{\d_1,\d_1\inv}$
and so rationally generated by $\set{\d_1}$
where $\d_1$ is the characteristic function of $\set1$. It is also standard that
the maximal ideal space of $A$ can be identified with the annulus
$\set {w\in\C\st\rho_-\le\ab w\le \rho_+}$. By standard results in dimension
theory $\O(A)$ has dimension 1 if $\rho_-=\rho_+$ and dimension 2 otherwise.
If $A$ has dense invertible group then so has $\hA$ and we must have
$\rho_-=\rho_+$ by Proposition 1.3, as $A$ is plainly symmetric.

Conversely suppose that $\rho_-=\rho_+$. Then the spectrum of $\d_1$ has
empty interior and Proposition 1.5 applies.\eop

\npar Recall that a {\it uniform algebra}, $A$, is a Banach subalgebra of $C(X)$
with the
uniform norm such that $A$ separates the points of $X$ (a compact Hausdorff space)
and contains the constants. A (unital) {\it Banach function algebra} satisfies the same axioms as
a uniform algebra except that the complete norm on the algebra need not be the uniform norm.
The Banach function algebra $A$ is called {\it natural} if
all of its characters are given by evaluation at points of $X$ and {\it
trivial} if $A=C(X)$.

In the following result we make use of the Arens-Royden theorem (see
for example $\ref{19}$,
p. 411) which
implies that whenever $A$ is a natural Banach function algebra
on a compact Hausdorff space, $X$, there is an isomorphism
$$G(A)\bigm/ e^A\to G(C(X))\bigm/ e^{C(X)}$$
and that these groups are isomorphic to the homotopy classes of maps $X\to S^1$.
These groups are also isomorphic to $H^1(X,\Z)$ the first \v Cech cohomology
group of $X$ with coefficients in $\Z$. (Further information about \v Cech cohomology groups
can be found in Section 10 of $\ref{23}$.) For example when $X$ is a
compact plane set then $H^1(X,\Z)=0$ if and only if $\C\setminus X$ is connected.

\thm 1.7. Let $A$ be a natural Banach function algebra
on a compact Hausdorff space, $X$, such
that $H^1(X,\Z)=0$. Then:
\item{(i)} If $X$ is locally connected and $A$ has dense invertible group then
$A$ is uniformly dense in $C(X)$ (and so trivial if $A$ is a uniform algebra).
\item{(ii)} If $X$ contains a locally connected closed subspace
whose topological dimension is at least two then $A$ does not have
a dense invertible group.

\pf We prove (ii); one can prove (i) by repeating parts of the argument.

Let $E$ be a non-empty, closed, locally connected subspace of $X$ of topological
dimension at least two and suppose for a contradiction that
$A$ has a dense invertible group.
By the comments above, $e^A=G(A)$. Let $B$
be the algebra of functions on $E$ which are restrictions of elements of $A$.
Then $B$ is isomorphic to the quotient of $A$ by the (closed) ideal of
functions which vanish on $E$. We regard $B$ as a Banach function algebra
on $E$ with respect to this quotient norm.

By hypothesis the set $\F=\set{f\in C(E)\st\hbox{ there exists }\tilde f\in e^A
\hbox{ with }\tilde f\vert_E=f}$ is dense in $B$. Now let $C$ be the closure
of $B$ in $C(E)$. It is easy to check that $C$ is a uniform algebra on $E$
such that $\F$ is dense in $C$. Clearly every element in $\F$ has a square root
in $C$ and
$E$ is locally connected so
it follows from \v Cirka's theorem (see
$\ref{23}$, p. 131-134) that $C=C(E)$.
Since $C(E)$ has a dense set of invertible elements we must
have $\dim E\le1$ by the fundamental result stated in the introduction,
a contradiction.\eop

\npar It is a famous open question (see Section 3) whether or not there exists
a non-trivial natural uniform algebra on the unit interval, $I$. In connection with this
we have the following corollary:

\cor 1.8. Let $A$ be a natural uniform algebra on the unit interval.
If $A$ has dense invertible group then $A=C(I)$.

\npar Now it follows quite generally and elementarily that elements of $A\setminus\bGA$ have
spectra whose interiors contain $0$. Thus the
set of space-filling curves of a non-trivial natural uniform algebra on the unit interval
must have interior in $A$. It has certainly been known for a long time ($\ref{14}$)
that every such uniform algebra must contain a function whose image has positive
Lebesgue measure.
One can also compare the result of \v Cirka ($\ref{15}$, p. 670)
which shows (using standard theory) that a non-trivial doubly
generated natural uniform algebra on $I$ must contain a dense set of functions 
whose images have non-empty interior.

Compare too the above with an unpublished result of Cole ($\ref{4}$)
which states that every non-trivial uniform algebra contains
a chain of subalgebras $B_1\subseteq B_2\subseteq\cdots$ such that for each
$n\in\N$, $B_n$ is isomorphic to the $n$-dimensional polydisc algebra.
We also
remind the reader of the well-known examples of (non-natural) uniform algebras $A$
on the unit interval such that every non-constant $f\in A$ has the interior of
$f(I)$ non-empty. (See for example $\ref{18}$, p. 200.)

The authors would like to thank to Brian Cole and John
Wermer for helpful discussions relating to the material in this section.

\eject


\vskip 20pt
\noindent{\bf 2. Integral Extensions}
\vskip 10pt
\noindent
As hinted at in the introduction denseness of invertible elements in
uniform algebras seems to be connected with the existence of roots of
monic polynomials over the algebra. Accordingly we mention the following result
stated by Grigoryan in $\ref{12}$: a uniform algebra $A$ on a compact space,
$X$, is trivial
if $C(X)$ is
finitely generated and is an integral extension of $A$.

Our main result is that integral extensions of commutative Banach algebras preserve the
property of having dense invertible group.

Recall (or see for example p. 369 of $\ref{19}$) that if $A$ is a commutative unital normed
algebra and $\a(x)$ is a monic polynomial over $A$ then the Arens-Hoffman extension, $\Aa$,
of $A$ by $\a(x)$ is $A[x]/\prin$ where $\prin$ denotes the principal ideal
in $A[x]$ generated by
$\a(x)$. There are infinitely many equivalent norms which make this algebra
an isometric extension of $A$. In fact, for a sufficiently
large fixed value of $t>0$, such an `Arens-Hoffman' norm is given by
$$\norm{\prin+\b(x)}=\sum_{j=0}^{n-1}\norm{b_j}t^j$$
where $\b(x)=\sum_{j=0}^{n-1}b_jx^j\ (b_0,\ldots,b_{n-1}\in A)$. (It can be shown,
as on p. 128 of $\ref{16}$ for example, that
every element of $\Aa$ has a unique representative of degree less than $n$ so this
is well defined.) We shall write $\xb$ for the coset $\prin+x$ from now on and
we shall assume that the norm on $\Aa$ is given by an Arens-Hoffman norm.

The proof of our main result relies on resultants; we now recall their definition.
Let $\a(x)=a_0+\cdots+
a_{n-1}x^{n-1}+x^n$ be a monic polynomial over a commutative ring $A$ and
$\b(x)=b_0+b_1x
\cdots+b_{n-1}x^{n-1}\in A[x]$.
By definition (see for example $\ref{16}$, p. 325)
$${\rm res}(\a(x), \b(x)):=\left\vert\matrix{
\overbrace{\matrix{1 &a_{n-1} &\cdots &a_2\cr
			0 &1 &\cdots &\cdots\cr
			\cdots  &\cdots &\cdots &\cdots\cr
			0  &\cdots &\cdots &1\cr}}^{n-1}
&\overbrace{\matrix{    a_1 &a_0 &0 &\cdots &0\cr
			a_2 &a_1 &a_0 &\cdots &0\cr
			\cdots &\cdots &\cdots &\cdots &\cdots \cr
			a_{n-1} &a_{n-2} &\cdots &\cdots &a_0 \cr}}^{n}	\cr
\matrix{b_{n-1} &b_{n-2} &\cdots &b_1\cr
	0 &b_{n-1} &\cdots &\cdots\cr
	\cdots &\cdots &\cdots &\cdots \cr
	0 &0 &\cdots &0\cr}
&\matrix{b_0 &0 &0 &\cdots &0\cr
	b_1 &b_0 &0 &\cdots &0\cr
	\cdots &\cdots &\cdots &\cdots &\cdots\cr
	b_{n-1} &b_{n-2} &\cdots &\cdots &b_0\cr}\cr}\right\vert$$
\vskip 5pt
It can be shown (see $\ref{2}$, A.IV \S 6.6.1) that $\prin+\b(x)$ is
invertible in $\Aa$ if and only if
$\res(\a(x),\b(x))$ is invertible in $A$.

We see from the above that, writing $c$ for $b_0$,
$\res(\a(x),\b(x))=P(c)=
p_0+p_1c+\cdots+p_{n-1}c^{n-1}+c^n$ for some $p_0,\ldots,p_{n-1}\in A$
which are polynomials in $b_1,\ldots,b_{n-1},a_0,\ldots,a_{n-1}$ only with
coefficients in $\C$.

We shall usually fix $\a(x)$ and allow $\b(x)$ to vary. We then
denote $\res(\a(x),\b(x))$ by $R_\a(\b(x))$ for $\b(x)\in A[x]$.

It is a standard fact ($\ref{11}$, p. 398)
that $R_\a(b_0+\cdots+b_{n-1}x^{n-1})$ is homogeneous of
degree $n$ in $b_0,\ldots,b_{n-1}$ and homogeneous of degree $n-1$ in $a_0,
\ldots,a_{n-1},1$.

Sometimes in the literature (for example $\ref{3}$ and $\ref{8}$) the integral
extensions considered are obtained by adjoining square roots. In this
case we have $\a(x)=x^2-a_0$ and then $R_\a(b_0+b_1x)=b_0^2-a_0b_1^2$.

\thm 2.1. Let $A$ be a commutative unital Banach algebra with dense invertible
group and $\a(x)$ a monic polynomial over $A$. Then $G(\Aa)$ is dense in $\Aa$.
\pf
We may assume $n\ge 2$. Let $\b(x)=b_0+b_1x+
\cdots+b_{n-1}x^{n-1}\in A[x]$ and $\e>0$. We have to show that there exists
$\tilde\b(x)\in A[x]$ with $\norm{\tilde\b(\xb)-\b(\xb)}<\e$ and $\tilde
\b(\xb)\in G(\Aa)$. In fact we shall show that by slightly perturbing
$b_0$ only, we can obtain a polynomial $\tilde\b(x)$ with $R_\a(\tilde\b(x))
\in G(A)$.

By the above comments, $P(c)=R_\a(c+b_1x+\cdots+b_{n-1}x^{n-1})$ is a polynomial
$p_0+\cdots+p_{n-1}c^{n-1}+c^n$ where $p_0,\ldots,p_{n-1}\in A$ are independent
of $c$. Consider the $n$ formal derivatives of $P$ as maps $A\to A$:
$$\eqalign{&P^{(0)}(c) = P(c);\cr
&P^\prime(c) = p_1+2p_2c+\cdots+nc^{n-1};\cr
&\qquad\vdots\cr
&P^{(n-1)}(c)=n!c.\cr}$$
Set $\tc 1=b_0$. Note that, trivially, $P^{(n-1)}$ is a local homeomorphism at
$\tc 1$.
Now let $1\le k<n$ and suppose that $\tc1,\ldots,\tc k\in A$ have been chosen so that
\item{(i)} $P^{(n-j)}$ is a local homeomorphism at $\tc j\ (1\le j\le k)$, and
\item{(ii)} $\norm{\tc j-\tcb {j-1}}<\e/n$ for $1<j\le k$.
We shall now show how to choose $\tcb{k+1}$ so that (i) and (ii) become
true with `$k$' replaced by `$k+1$'.

It is easy to see from the inverse function theorem for Banach spaces
(see for example $\ref{1}$, Chapter 7 Theorem 8) that $P^{(n-(k+1))}$ is a local
homeomorphism at $a\in A$ if $P^{(n-k)}(a)$ is invertible. (This fact is also stated
in $\ref 7$.)

By hypothesis, $P^{(n-k)}$ is a local homeomorphism at $\tc k$ so some open
neighbourhood of $\tc k$ is mapped onto an open set in $A$. Since $G(A)$ is
dense in $A$ there is some $\tcb{k+1}\in B_A(\tc k,\e/n)$ with
$P^{(n-k)}(\tcb{k+1})\in G(A)$.

Thus $\tcb{k+1}$ has the required
properties and by induction we can choose $\tilde c_{n-1},\ldots,\tilde c_0\in A$
with $\norm{\tilde c_k-\tilde c_{k-1}}<\e/n\ (k=1,\ldots,n-1)$, $\tc 1=b_0$,
and $P$ a local homeomorphism at $\tilde c_0$.

Again, since $P$ is a local homeomorphism at $\tilde c_0$, we can find $\tilde b_0
\in B_A(\tilde c_0,\e/n)$ with $P(\tilde b_0)\in G(A)$. Since
$$\norm{\tilde b_0-b_0}\le\norm{\tilde b_0-\tilde c_0}+
\norm{\tilde c_0-\tilde c_1}+\cdots+\norm{\tilde c_{n-2}-\tilde c_{n-1}}<\e,$$
the result is proved.\eop

\npar
We also remark that the method of proving Theorem 2.1 gives another way to see
that if $A$ is a commutative unital Banach algebra with $\bGA=A$ then for
every $n\in\N$, $M_n(A)$ also has dense invertible group. In fact we can approximate
any $B=[b_{ij}]\in M_n(A)$ by an invertible matrix by perturbing only $n$ entries
$b_{i_1,j_1},\ldots, b_{i_n,j_n}$ provided that $k\mapsto i_k,j_k$ are both
permutations. We leave the details to the reader.

\cor 2.2. Let $A$ and $B$ be commutative
normed algebras and suppose that $B$ is an integral
extension of $A$. Suppose that $A$ is a full subalgebra of its completion, $\tA$,
and that $A$ has dense invertible group. Then $B$ has dense invertible group.

\pf It is sufficient to prove the case when $B$ is an Arens-Hoffman extension
of $A$ for if $C$ is a normed integral extension of $A$ and $c\in A$ is a root
of the monic polynomial $\a(x)\in A[x]$ then there is a continuous unital
homomorphism $\theta\;\colon\,\Aa\to C$ with $\theta(\xb)=c$; the result quickly
follows from this.

Let $b\in \Aa$ and let $\e>0$.
It is not hard to show that the universal property of Arens-Hoffman extensions which
has just been mentioned allows us to identify
$(\tilde A)_\a$ with $\widetilde{\Aa}$. By Theorem 2.1
there exists $b\pri\in G((\tA)_\a)$ with $\norm{\bpr-b}<\e/2$.
Suppose $\bpr=\sum_{k=0}^{n-1}\bpr_k\xb^k$ where $n$ is the degree of $\a(x)$ and
$\bpr_0,\ldots,\bpr_{n-1}\in \tA$. Since $G((\tA)_\a)$ is open (the Arens-Hoffman
extension of a Banach algebra is complete) we can perturb $\bpr_0,\ldots,\bpr_{n-1}$
to obtain $\bppr_0,\ldots,\bppr_{n-1}\in A$ so that $\bppr=\sum_{k=0}^{n-1}{\bpr}^\prime_k\xb^k$
is invertible in $(\tA)_\a$ and $\norm{\bppr-\bpr}<\e/2$. But now
$R_\a(\bppr)\in G(\tA)\cap A= G(A)$ so $\bppr$ is invertible in $\Aa$
and $\norm{b^{\prime\prime}-b}<\e$. \eop

\npar The following theorem follows directly from this corollary:

\thm 2.3. Let $B$ be a commutative Banach algebra which is an
integral extension of the commutative Banach algebra $A$. If
$A$ has dense invertible group then so has $B$.

\npar The converse of Theorem 2.1 seems harder to investigate. However the method of
the proof gives at least partial information in this direction:

\prop 2.4. Suppose that $A$ is a commutative unital Banach algebra and $\a(x)$
is a monic polynomial of degree $n$ over $A$. If $\Aa$ has dense invertible group
then $\set{b\in A\st\hbox{there exists } a\in A\hbox{ such that } b=a^n}$, the set
of $n$th powers, is contained in the closure of $G(A)$.

\pf Fix $a\in A$ and let $\e>0$. Let the norm parameter for the Arens-Hoffman
extension be $t>0$. Then there exists $\b_\e(x)
=b_{\e,0}+\cdots+b_{\e,n-1}x^{n-1}\in A[x]$ such $R_\a(\b_\e(x))\in G(A)$
and
$$\norm{\b_\e(x)-a}=\norm{b_{\e,0}-a}+\sum_{j=1}^{n-1}\norm{b_{\e,j}}t^j<\e.$$
Now as we mentioned before, $R_\a(b_0+\cdots+
b_{n-1}x^{n-1})$ is homogeneous of degree $n$ in $b_0,
\ldots,b_{n-1}$. Therefore, writing $R_\a(\b(x))=P(b_0)=\sum_{j=0}^{n-1}p_jb_0^j
+b_0^n$ as in Theorem 2.1, we have that each coefficient $p_j\ (j=0,\ldots,
n-1)$ is a sum of elements of $A$ each of which has $b_k$ as a factor
for some $k\in\set{1,\ldots,n-1}$.

Thus, letting $\e\to0$, we obtain invertible elements $R_\a(\b_\e(x))$ in $A$ which
tend to $a^n$.\eop

\npar It is clear that the property of having a dense invertible group passes to
quotients and completions of normed algebras. The following is also clear.

\prop 2.5. Let $A$ be a direct limit of normed algebras with dense invertible groups
and where the connecting homomorphisms are unital isometric monomorphisms. Then
$A$ has dense invertible group.

\pf This is clear.\eop

\npar We mention the consequences of these results
for some of the examples of uniform algebras constructed by Karahanjan, Cole, and Feinstein
($\ref{17}$, $\ref{3}$, $\ref{8}$).
Specifically we are referring to those obtained
by taking Banach algebra direct limits of systems
of uniform algebras in which the intermediate algebras are completions
(or quotients followed by completions) of algebraic extensions of their
predecessors. In many of these the initial algebra has dense invertible group
and so the final algebra also has this property. In particular
our results show that
there is no need to invoke the theory of `dense thin
systems' as developed in $\ref{17}$.


\vskip 20pt
\noindent{\bf 3. Open Questions}
\vskip 10pt

\noindent{\bf 1.} The question of whether or not there is more than one natural uniform
algebra on the unit interval is now about fifty years old. It seems to
have been first formally asked by Gelfand in 1957 ($\ref{10}$). The answer is not
known even if we assume that the algebra is regular. See $\ref{9}$ and the
references cited there for material on this problem.

\npar{\bf 2.} Let $A$ be a commutative Banach algebra with
maximal ideal space $\O$ and \v Silov boundary $S$.
A lack of examples describing $\bGA$ leads to the following
conjectures:
\item{(a)} If $\bGA=A$ then $S=\O$.
\item{(b)} If $\bGA=A$ then $C(\O)$ has dense invertible group.
\item{(c)} If $S=\O$ and the invertible elements of $C(\O)$ are dense in $C(\O)$
then $\bGA= A$.

It follows from standard properties about the maximal ideal space
and \v Silov boundary of the uniform closure of a Banach function algebra
(see for example $\ref6$, p. 447) that we may assume $A$ is a uniform algebra in
Questions 2(a) and 2(b). The answer to 2(c) may depend on the category.

\npar{\bf 3.} The converse of Theorem 2.1 remains open.

\npar{\bf 4.} We do not know if the word `full' is redundant in Corollary 2.2.

\vfill\eject


\vskip 20pt
\noindent{\bf References }

\vskip 10pt
\noindent$\ref{1}$\ Bollob\'as, B. (1999): `Linear Analysis', 2nd edition, Cambridge:
Cambridge University Press.

\vskip 10pt
\noindent$\ref{2}$\ Bourbaki, N. (1988) `Algebra II, Chapters 4-7', US: Springer-Verlag.

\vskip 10pt
\noindent$\ref{3}$\ Cole, B. J. (1968): `One-Point Parts and the Peak Point
Conjecture', PhD Thesis, Yale University.

\vskip 10pt
\noindent$\ref{4}$\ Cole, B. J. (2002): {\sl Private communication.}

\vskip 10pt
\noindent$\ref{5}$\ Corach, G. and Su\' arez, F. D. (1988): 
`Thin Spectra and Stable Range Conditions', {\sl J. Funct. Anal.}, 81, 432-442.

\vskip 10pt
\noindent$\ref{6}$\ Dales, H. G. (2001): `Banach Algebras and Automatic
Continuity', New York: Oxford University Press Inc.

\vskip 10pt
\noindent$\ref7$\ Falc\'on-Rodr\'\i guez, C. M. (1988) `Sobre la densidad del grupo de
los elementos invertibles de un \'algebra uniforme', {\sl Revista Ciencias
Matem\'aticas}, IX, no. 2, 11-17.

\vskip 10pt
\noindent$\ref{8}$\ Feinstein, J. F. (1992): `A Non-Trivial, Strongly Regular Uniform
Algebra', {\sl J. Lond. Math. Soc.}, 45, no. 2, 288-300.

\vskip 10pt
\noindent$\ref{9}$\ Feinstein, J. F. and Somerset, D. W. B. (1999):
`Strong Regularity for Uniform Algebras', {\sl Contemp. Math.}, 232, 139-149.

\vskip 10pt
\noindent$\ref{10}$\ Gelfand, I. M. (1957): `On the Subrings of a Ring of
Continuous Functions', {\sl Usp. Mat. Nauk}, 12, no. 1, 249-251.

\vskip 10pt
\noindent$\ref{11}$\ Gelfand, I. M., Kapranov, M. M., and Zelevinsky, A. V. (1994)
`Discriminants, Resultants, and Multidimensional Determinants', United
States of America: Birkh\"auser Boston.

\vskip 10pt
\noindent$\ref{12}$\ Grigoryan, S. A. (1984): `Polynomial Extensions of Commutative
Banach Algebras', {\sl Russian Math. Surveys}, 39, no. 1, 161-162.

\vskip 10pt
\noindent$\ref{13}$\ Hatori, O. and Miura, T.: (1999) `On a Characterization of
the Maximal Ideal Spaces of Commutative $C^*$-Algebras in
which Every Element is the Square of Another', {\sl Proc. Am. Math. Soc.}, 128,
no. 4, 1185-1189.

\vskip 10pt
\noindent$\ref{14}$\ Helson, H. and Quigley, F. (1957): `Existence of Maximal Ideals
in Algebras of Continuous Functions', {\sl Proc. Am. Math. Soc.}, 8, 115-119.

\vskip 10pt
\noindent$\ref{15}$\ Henkin, G. M. and \v Cirka, E. M. (1976):
`Boundary Properties of Holomorphic Functions of Several Complex
Variables', {\sl J. Soviet Math.}, 5, no. 5, 612-687.

\vskip 10pt
\noindent$\ref{16}$\ Jacobson, N. (1996) `Basic Algebra I' (2nd ed.) New York:
W. H. Freeman and Company.

\vskip 10pt
\noindent$\ref{17}$\ Karahanjan, M. I. (1979): `Some algebraic characterizations of the
algebra of all continuous functions on a locally connected compactum',
{\sl Math. USSR Sb.}, 35, 681-696.

\vskip 10pt
\noindent$\ref{18}$\ Leibowitz, G. M. (1970) `Lectures on Complex Function Algebras',
United States of America: Scott, Foresman and Company.

\vskip 10pt
\noindent$\ref{19}$\ Palmer, T. W. (1994) `Banach Algebras and the General
Theory of *-Algebras' (vol. 1), Cambridge: CUP.

\vskip 10pt
\noindent$\ref{20}$\ Pears, A. R. (1975) `Dimension Theory of General Spaces',
Cambridge: CUP.

\vskip 10pt
\noindent$\ref{21}$\ Rieffel M. A. (1983): `Dimension and Stable Rank in $K$-Theory of
$C^*$-Algebras', {\sl Proc. Lond. Math. Soc.}, 46, no. 3, 577-600.

\vskip 10pt
\noindent$\ref{22}$\ Robertson, G. (1976): `On the Density of the Invertible Group
in $C^*$-Algebras', {\sl Proc. Edinb. Math. Soc.}, 20, 153-157.

\vskip 10pt
\noindent$\ref{23}$\ Stout, E. L. (1973) `The Theory of Uniform Algebras',
Tarrytown-on-Hudson, New York: Bogden and Quigley Inc.

\vskip 40pt

\npar
T. W. Dawson

\vskip 4pt

\noindent E-mail: {\tt pmxtwd@nottingham.ac.uk}

\vskip 14pt

\par\noindent Dr. J. F. Feinstein

\vskip 4pt

\noindent E-mail: {\tt Joel.Feinstein@nottingham.ac.uk}

\npar
Division of Pure Mathematics,
\par\noindent School of Mathematical Sciences,
\par\noindent University of Nottingham,
\par\noindent University Park,
\par\noindent Nottingham,
\par\noindent NG7 2RD,
\par\noindent UK.

\end